\documentclass[12pt]{article}

\usepackage{listings}
\usepackage{url}

\usepackage[english]{babel}
\usepackage{setspace}
\usepackage[letterpaper,top=2cm,bottom=2cm,left=3cm,right=3cm,marginparwidth=1.75cm]{geometry}
\usepackage[dvipsnames,table,xcdraw]{xcolor}

\newcommand{\bigconcchoose}[1]{\def\bigconcsize{}%
  \ifx#1\displaystyle
    \let\bigconcsize\Big
  \else
    \ifx#1\textstyle
      \let\bigconcsize\big
    \fi
  \fi#1}

\providecommand{\keywords}[1]
{
  \small	
  \textbf{\textit{Keywords---}} #1
}

\usepackage{amsmath}
\usepackage{graphicx}
\usepackage[colorlinks=true, allcolors=blue]{hyperref}
\usepackage{amsfonts}
\usepackage{amsthm}
\usepackage{amssymb}

\newtheorem{lemma}{Lemma}
\newtheorem{theorem}[lemma]{Theorem}
\newtheorem{proposition}[lemma]{Proposition}
\newtheorem{corollary}[lemma]{Corollary}
\newtheorem{definition}[lemma]{Definition}

\usepackage{graphicx}
\usepackage{qtree}
\usepackage{amsmath}
\usepackage{amsthm}
\usepackage{amssymb}
\usepackage[dvipsnames]{xcolor}
\usepackage{tikz,tkz-euclide}
\usetikzlibrary{decorations.pathreplacing,angles,quotes}
\usetikzlibrary{calc}
\usetikzlibrary{snakes}
\usetikzlibrary{arrows}
\usetikzlibrary{arrows.meta}
\usetikzlibrary{positioning}
\usepackage{mathtools}
\usepackage{algpseudocode}
\usepackage{mathrsfs}
\usepackage{xpatch}
\usetikzlibrary{nfold}
\usepackage[font={footnotesize,stretch=1.3}]{caption}

\DeclareMathOperator{\diam}{diam}

\DeclareMathOperator{\abs}{abs}
\DeclareMathOperator{\True}{True}
\DeclareMathOperator{\False}{False}

\usepackage{hyperref} 
\hypersetup{
    colorlinks,
    citecolor=blue,
    filecolor=black,
    linkcolor=blue,
    urlcolor=blue
}

\title{Diameter Constraints in 2-distance Graphs}

\author{Oleksiy Al-saadi$^{\dag}$, Joseph Natal$^{\ddagger}$  \\
        \small $^{\dag}$Sonoma State University, {\nolinkurl{alsaadio@sonoma.edu}} \\
        \small $^{\ddagger}$Karlsruhe Institute of Technology, {\nolinkurl{joseph.natal@kit.edu}} \\
}

\begin{document}

\maketitle

\begin{abstract}
For any finite, simple graph $G = (V,E)$, its $2$-distance graph $G_2$ is a graph having the same vertex set $V$ where two vertices are adjacent if and only if their distance is $2$ in $G$. Connectivity and diameter properties of these graphs have been well studied. For example, it has been shown that if $\diam(G) = k \geq 3$ then $\lceil \frac{1}{2} k \rceil \leq \diam(G_2)$, and that this bound is sharp. In this paper, we prove that $\diam(G_2) = \infty$ (that is, $G_2$ is disconnected) or otherwise $\diam(G_2) \leq k + 2$. In addition, we show that this inequality is sharp for any even $k$, a result that we verify for some higher orders through judicious use of a \textsc{sat} solver.
\end{abstract}




\section{Introduction}  

Given a finite, simple graph $G$, we let $V$ and $E$ denote the vertex set and the edge set of $G$, respectively. By $d_G(u,v)$ we denote \emph{geodesic distance}, the length of a shortest path, between the vertices $u$ and $v$ in the graph $G$. The \emph{diameter} of $G$, also denoted as $\diam(G)$, is the greatest geodesic distance between any two vertices in $G$. If $G$ is disconnected, then $\diam(G) = \infty$. The \emph{$k$-distance operator} $D_k$ takes as input a simple graph $G$ and outputs its \emph{$k$-distance graph}, a graph having the same vertex set $V$ but where any pair $(u,v)$ has $d_{G_k}(u,v) = 1$ if and only if $d_G(u,v) = k$. Thus, $D_k(G)$ is the $k$-distance graph of $G$. A graph $Y$ is \emph{$k$-distance} if there exists a graph $X$ such that $Y \cong D_k(X)$ (meaning that $Y$ is isomorphic to $D_k(X)$). When ambiguity is impossible, we simply write $G_k$ to refer to $D_k(G)$.

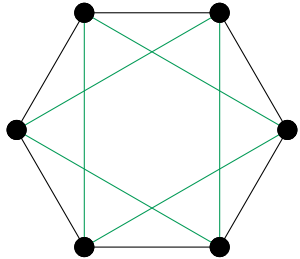
\begin{figure}[h]
    \centering
        \begin{tikzpicture}[scale=.9, every node/.style={scale=.6}]

        \tikzstyle{vertex}=[circle, fill=white]
        \begin{scope}\LARGE
        \end{scope}
        
        \tikzstyle{vertex}=[circle, fill=black]
        \tkzDefPoint(0,0){O}\tkzDefPoint(2,0){A}
        \tkzDefPointsBy[rotation=center O angle 360/6](A,B,C,D,E){B,C,D,E,F}
        \tkzDrawPoints[fill =black,size=12,color=black](A,B,C,D,E,F)
        

        \tikzset{decoration={snake,amplitude=.3mm,segment length=3mm, post length=0mm,pre length=0mm}}

        \draw(A)--(B)--(C)--(D)--(E)--(F)--(A);
        \draw[ForestGreen](A)--(C)--(E)--(A);
        \draw[ForestGreen](B)--(D)--(F)--(B);

        \tikzstyle{vertex}=[circle, fill=black]
        \tkzDefPoint(0,0){O}\tkzDefPoint(2,0){A}
        \tkzDefPointsBy[rotation=center O angle 360/6](A,B,C,D,E){B,C,D,E,F}
        \tkzDrawPoints[fill =black,size=12,color=black](A,B,C,D,E,F)

        \end{tikzpicture}
        \caption{Let $G = (V,E)$ be the graph depicted by the black vertices and edges. Then, let the green edges represent $E(G_2)$. Observe that $G$ is an even cycle and that $G_2$ contains two disconnected triangles.}
        \label{fig:evenDisc}

\end{figure}

The general class of $k$-distance graphs was first introduced by Harary et al. \cite{Harary} and are also known under the name of \emph{exact distance powers} \cite{Bai,Foucaud}. Special cases for fixed $k$ have been studied for their connectivity, diameter, and periodicity properties \cite{Connectivity2Dist,Jafari,Prisner,BohdanZelinka2000,Erick}. Known results in the literature are numerous. For example, the cubic self $2$-distance graphs do not exist \cite{NoCubic}. On the other hand, Azimi et al. \cite{simplePathCycle} showed the complete set of graphs whose $2$-distance graphs are simple paths or cycles. Graphs for which $G \cong G_2$ are particularly interesting. As a simple example, it is easily verified that any odd cycle has an isomorphic $2$-distance graph---whereas the $2$-distance graph of an even cycle is disconnected (See Fig. \ref{fig:evenDisc}). General characterizations were given by Ching and Garces \cite{ThreeChar} whereas Gaar and Krenn \cite{GAAR2023181} characterized regular $2$-distance graphs. Recognition under special constraints was shown to be polynomial-time by Bai et al. \cite{Bai}. 

Recently, Jafari and Musawi \cite{Jafari} conjectured that if $\diam(G) = k \geq 3$ then $\diam(G_2) \leq k + 2$. Under the standard assumption that $\diam(G) = \infty$ when $G$ is disconnected, we will prove that the statement is true when $\diam(G_2) \neq \infty$.

\section{Preliminaries}

Our graph theory notation basically follows that of Golumbic \cite{Golumbic}. By $G$ we denote a simple, undirected, and finite graph. The vertex set of $G$ is denoted by $V$ and the edge set of $G$ is denoted by $E$. We write $G[S]$ to refer to the graph induced by the vertex set of $S$. For a graph $H$, we let $V(H)$ and $E(H)$ denote the vertex set and the edge set of $H$, respectively. For any $x \in V$, let $N_G(x) = \{y : \{x,y\} \in E\}$ be the (open) neighborhood of $x$ and $N_G[x] = N(x) \cup \{x\}$ the closed neighborhood of $x$. We write $x \in N_G(S)$ to say that $x$ is adjacent to some vertex in $S$ while $x \not\in S$. We write $x \in N_G[S]$ if $x \in N_G(S) \cup S$. We denote by $C_{N_G}(S)$ the \emph{common neighborhood} of $S$ in $G$. i.e.
\[
C_{N_G}(S) = \bigcap_{v \in S} N_G(v).
\]
A sequence of vertices $P = \langle u = x_0 x_1 \cdots x_{k-1} x_k = v \rangle$ is called a \emph{walk} if $x_i \in N_G[x_{i+1}]$ for all $i \in \{0,1,\dots,k-1\}$. The \emph{length} of a walk is the number of edges it contains. We say $P$ is a \emph{path} if the vertices $x_0,x_1, \dots ,x_k$ are all distinct. A walk with endpoints $u$ and $v$ may be called a $u,v$-walk. If $v \in N_G(v')$, then we denote by $P \mathord{-} v'$ the walk $\langle u = x_0 x_1 \cdots x_{k-1}  x_k = v v' \rangle$. Similarly, if $P = \langle x_0 x_1  \cdots x_k \rangle$ and $P' = \langle y_0 y_1 \cdots y_q \rangle$ are walks and $x_k \in N_G(y_0)$, then we let $P \mathord{-} P' = \langle x_0 \cdots x_k y_0 \cdots y_q \rangle$. We write $P[x_i,x_{j}]$ where $(0 \leq i \leq j \leq k)$ for the subwalk $\langle x_i  \cdots x_j \rangle$ of $P$. It is well known that we can extract a $u,v$-path from the subset of every $u,v$-walk. A \emph{triangle} is a set of three vertices such that each vertex is adjacent to the other two.

\section{Diameter Bounds}

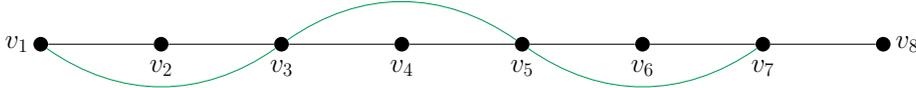
\begin{figure}[h]
    \centering
        \begin{tikzpicture}[scale=.8, every node/.style={scale=.5}]

        \tikzstyle{vertex}=[circle, fill=white]
        \begin{scope}\LARGE
        \node[vertex] (a1) at (-.4,0) {$v_1$};
        \node[vertex] (a2) at (2,-.4) {$v_2$};
        \node[vertex] (a3) at (4,-.4) {$v_3$};
        \node[vertex] (a4) at (6,-.4) {$v_4$};
        \node[vertex] (a5) at (8,-.4) {$v_5$};
        \node[vertex] (a6) at (10,-.4) {$v_6$};
        \node[vertex] (a7) at (12,-.4) {$v_7$};
        \node[vertex] (a8) at (14.4,0) {$v_8$};

        \end{scope}
        
        \tikzstyle{vertex}=[circle, fill=black]
        
        \node[vertex] (a1) at (0,0) {$ $};
        \node[vertex] (a2) at (2,0) {$ $};
        \node[vertex] (a3) at (4,0) {$ $};
        \node[vertex] (a4) at (6,0) {$ $};
        \node[vertex] (a5) at (8,0) {$ $};
        \node[vertex] (a6) at (10,0) {$ $};
        \node[vertex] (a7) at (12,0) {$ $};
        \node[vertex] (a8) at (14,0) {$ $};
        
        \tikzset{decoration={snake,amplitude=.3mm,segment length=3mm, post length=0mm,pre length=0mm}}
        
        \draw[black](a1)--(a2)--(a3)--(a4)--(a5)--(a6)--(a7)--(a8);

        \draw[ForestGreen,bend right=35](a1)edge(a3);
        \draw[ForestGreen,bend right=-35](a3)edge(a5);
        \draw[ForestGreen,bend right=35](a5)edge(a7);
            
        \end{tikzpicture}
        \caption{Let $P = \langle v_1 v_2 \cdots v_8 \rangle$ in $G$, shown with black edges. The green edges are the edges of $D_2(P)$ by Definition \ref{def:D2PHalved}.}
        \label{fig:halved}

\end{figure}

In this section, we show that the diameter of $G$ imposes strong constraints on the diameter of $G_2$ (see Theorem \ref{lem:k+2Upper}). We need an easy way to refer to paths between $G$ and its $2$-distance graph in our proofs. For a basic example of the following definition, see Fig. \ref{fig:halved}.

\begin{definition}
Let $P = \langle v_1 v_2 \cdots v_{k-1} v_k \rangle$ be a walk in $G$ and let $p = k - (k \mod 2)$. If for each odd $i \in [1,k-2]$ we have that $P[v_i,v_{i+2}]$ is induced, then there exists a walk $\langle v_1 v_3 \cdots v_{p-2} v_{p} \rangle $ in $G_2$ that we denote by $D_2(P)$. 
\label{def:D2PHalved}
\end{definition}

We will prove our main theorem by showing that a graph $G_2$ having an exceedingly high diameter causes $G$ to contain the complement of a path. This high edge-density subgraph will lead to a contradiction to $\diam(G)$. This result builds on a sequence of lemmas which we will now begin.

\begin{lemma}
Let $G$ be a graph with $\diam(G) = k \geq 3$, let $P_2 = \langle a_1 a_2 \cdots a_{k+3} a_{k+4} \rangle$ be a shortest path in $G_2$, and let $\ell = \lceil \frac{k+4}{2} \rceil$. Then, $a_1,a_{k+4} \in N_G(a_{\ell +1})$.
\label{lem:upperBound2Dist}
\end{lemma}

\begin{proof}
For $P_2$ to exist in $G_2$, the graph $G$ must contain a walk taking the form of $P = \langle a_1 b_1 a_2 b_2 a_3 \cdots a_{k+3} b_{k+3} a_{k+4} \rangle$. We will first show that the shortest paths between vertices of $V(P_2)$ in $G$ satisfy special constraints.

No shortest path in $G$ has length exceeding $\diam(G)$. In particular, the subwalks $P[a_1,a_{\ell+1}]$ and $P[a_{\ell+1},a_{k+4}]$ are not shortest paths, otherwise their length exceeds $\diam(G) = k$. Thus, there exist shortest paths $R = \langle a_1 = r_1 r_2 \cdots r_p = a_{\ell+1} \rangle$ and $R' = \langle a_{\ell+1} = r_1' r_2' \cdots r'_{p'} = a_{k+4} \rangle$ in $G$. In the remainder of the proof, we will frequently recall the fact that the length of $D_2(R)$ is at most $\lfloor \frac{k}{2} \rfloor$, which arises from the understanding that $D_2(R)$ is at most half the length of $R$ by Definition \ref{def:D2PHalved}.

We claim that $R$ and $R'$ have an odd length. Suppose instead that $L = R$ (\emph{resp.} $L = R'$) has an even length. Note that $\mu = (\ell + 1) - 1$ is the difference between the indices of $a_{\ell+1}$ and $a_1$. Moreover, $\omega = k+4 - (\ell+1)$ is the difference between the indices of $a_{k+4}$ and $a_{\ell+1}$. We have that $D_2(L)$ is shorter than both $P_2[a_1,a_{\ell+1}]$ and $P_2[a_{\ell+1},a_{k+4}]$ because $\lfloor \frac{k}{2} \rfloor$ (the greatest possible length of $D_2(L)$) is strictly less than both $\mu$ and $\omega$, respectively, a fact easily verified for any value of $k \geq 3$. This is a contradiction to the assumption that $P_2$ is a shortest path. Fig. \ref{fig:NoEvenPaths} shows an example for $k = 4$ and $L = R'$. We have proven our claim.

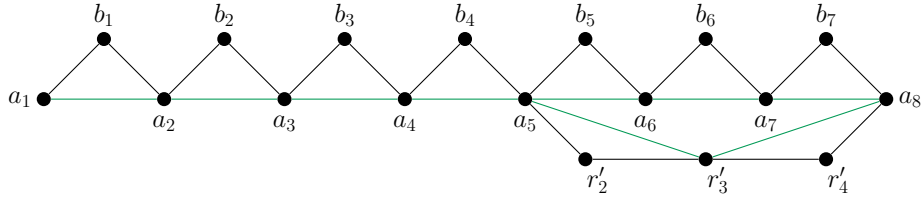
\begin{figure}[h]
    \centering
        \begin{tikzpicture}[scale=.8, every node/.style={scale=.48}]

        \tikzstyle{vertex}=[circle, fill=white]
        \begin{scope}\LARGE
        \node[vertex] (a1) at (-.4,0) {$a_1$};
        \node[vertex] (b1) at (1,1.4) {$b_1$};
        \node[vertex] (a2) at (2,-.4) {$a_2$};
        \node[vertex] (b2) at (3,1.4) {$b_2$};
        \node[vertex] (a3) at (4,-.4) {$a_3$};
        \node[vertex] (b3) at (5,1.4) {$b_3$};
        \node[vertex] (a4) at (6,-.4) {$a_4$};
        \node[vertex] (b4) at (7,1.4) {$b_4$};
        \node[vertex] (a5) at (8,-.4) {$a_5$};
        \node[vertex] (b5) at (9,1.4) {$b_5$};
        \node[vertex] (a6) at (10,-.4) {$a_6$};
        \node[vertex] (b6) at (11,1.4) {$b_6$};
        \node[vertex] (a7) at (12,-.4) {$a_7$};
        \node[vertex] (b7) at (13,1.4) {$b_7$};
        \node[vertex] (a8) at (14.4,0) {$a_8$};

        \node[vertex] (r1) at (9.2,-1.35) {$r'_2$};
        \node[vertex] (r2) at (11.2,-1.35) {$r'_3$};
        \node[vertex] (r3) at (13.2,-1.35) {$r'_4$};
        \end{scope}
        
        \tikzstyle{vertex}=[circle, fill=black]
        
        \node[vertex] (a1) at (0,0) {$ $};
        \node[vertex] (b1) at (1,1) {$ $};
        \node[vertex] (a2) at (2,0) {$ $};
        \node[vertex] (b2) at (3,1) {$ $};
        \node[vertex] (a3) at (4,0) {$ $};
        \node[vertex] (b3) at (5,1) {$ $};
        \node[vertex] (a4) at (6,0) {$ $};
        \node[vertex] (b4) at (7,1) {$ $};
        \node[vertex] (a5) at (8,0) {$ $};
        \node[vertex] (b5) at (9,1) {$ $};
        \node[vertex] (a6) at (10,0) {$ $};
        \node[vertex] (b6) at (11,1) {$ $};
        \node[vertex] (a7) at (12,0) {$ $};
        \node[vertex] (b7) at (13,1) {$ $};
        \node[vertex] (a8) at (14,0) {$ $};

        \node[vertex] (r1) at (9,-1) {$ $};
        \node[vertex] (r2) at (11,-1) {$ $};
        \node[vertex] (r3) at (13,-1) {$ $};

        \tikzset{decoration={snake,amplitude=.3mm,segment length=3mm, post length=0mm,pre length=0mm}}
        
        \draw(a1)--(b1)--(a2)--(b2)--(a3)--(b3)--(a4)--(b4)--(a5)--(b5)--(a6)--(b6)--(a7)--(b7)--(a8);

        \draw[ForestGreen](a1)--(a2)--(a3)--(a4)--(a5)--(a6)--(a7)--(a8);

        \draw[ForestGreen](a5)--(r2)--(a8);

        \draw(a5)--(r1)--(r2)--(r3)--(a8);

            
        \end{tikzpicture}
        \caption{Let $k = 4$. So, $P_2$ has length $7$ and $\ell+1 = 5$. Let the black edges denote $E(G)$ and let the green edges denote $E(G_2)$. A shortest $a_5,a_8$-path of length $4$ in $G$ causes $d_{G_2}(a_5,a_8) = 2$, a contradiction.}
        \label{fig:NoEvenPaths}

\end{figure} 

There is an important distinction between odd and even $k$. When $k$ is odd, then $R$ has length at most $k$. When $k$ is even, then $R$ has length at most $k - 1$ because we have shown that $R$ may not have even length. Let $k' = k$ when $k$ is odd, but $k' = k - 1$ when $k$ is even. Later in the proof, we will need to use $k'$ in order to reach a contradiction.

If $R$ and $R'$ both have length $1$, then we are done with the proof. Thus, at least one of these paths has length $3$ or greater. We will complete this proof by showing that a contradiction arises for any length of $R$. Note that $r_{p-1}$ is the last endpoint of $D_2(R)$ because $R$ has odd length. Consider two cases based on the incidence of the second to last vertex in $R$ to $R'$.

\textbf{Case 1:} $r_{p-1}$ is in the closed neighborhood of every vertex in $R'$ in $G$. If $a_{k+4} = r_{p-1}$ then $D_2(R)$ is an $a_1,a_{k+4}$-path of length at most $\lfloor \frac{k}{2} \rfloor$, a contradiction to the assumption that $d_{G_2}(a_1,a_{k+4}) = k + 3$. Thus, $a_{k+4} \neq r_{p-1}$. Next, if $X = \langle a_{k+4} r_{p-1} a_{\ell+1} \rangle$ is an induced path in $G$ then $D_2(X) = \langle a_{k+4} a_{\ell+1} \rangle$ gives $d_{G_2}(a_{k+4},a_{\ell+1}) = 1$, a contradiction. So, $X$ is not induced, meaning that $X$ is a triangle in $G$, implying $a_{k+4} \in N_G(a_{\ell+1})$. Now $\langle a_{\ell+1} a_{k+4} \rangle$ is a shortest path in $G$, and therefore $R'$ has length $1$ (i.e. $R' = \langle a_{\ell+1} a_{k+4} \rangle$). Since $R$ and $R'$ cannot both have length $1$, and neither has even length, it must be that $R$ has length $3$ or greater.

Let $B = C_{N_G}(a_{\ell},a_{\ell+1})$. For any $b \in B$, suppose that $b \not\in N_G[r_{p-1}] \cup N_G[a_{k+4}]$. Observe that $D_2(R \mathord{-} b \mathord{-} R')$ has length at most $\lfloor \frac{k}{2} \rfloor + 2$, a contradiction to the assumption that $d_{G_2}(a_1,a_{k+4}) = k + 3$. Fig. \ref{fig:Case1_XandB} demonstrates this contradiction where $k = 4$ and $b = b_{\ell}$. Thus, $b \in N_G[r_{p-1}] \cup N_G[a_{k+4}]$ for any $b \in B$.

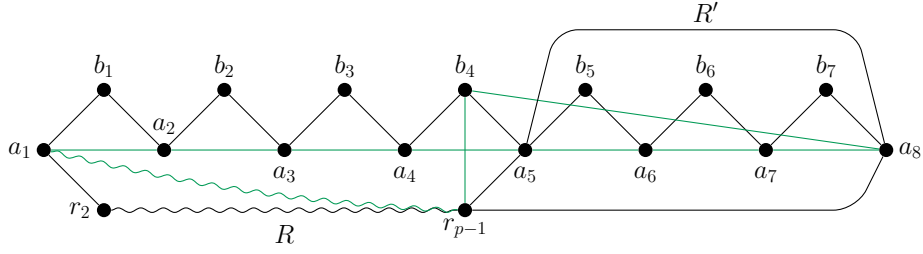
\begin{figure}[h]
    \centering
        \begin{tikzpicture}[scale=.8, every node/.style={scale=.48}]

        \tikzstyle{vertex}=[circle, fill=white]
        \begin{scope}\LARGE
        \node[vertex] (a1) at (-.4,0) {$a_1$};
        \node[vertex] (b1) at (1,1.4) {$b_1$};
        \node[vertex] (a2) at (2,.4) {$a_2$};
        \node[vertex] (b2) at (3,1.4) {$b_2$};
        \node[vertex] (a3) at (4,-.4) {$a_3$};
        \node[vertex] (b3) at (5,1.4) {$b_3$};
        \node[vertex] (a4) at (6,-.4) {$a_4$};
        \node[vertex] (b4) at (7,1.4) {$b_4$};
        \node[vertex] (a5) at (8,-.4) {$a_5$};
        \node[vertex] (b5) at (9,1.4) {$b_5$};
        \node[vertex] (a6) at (10,-.4) {$a_6$};
        \node[vertex] (b6) at (11,1.4) {$b_6$};
        \node[vertex] (a7) at (12,-.4) {$a_7$};
        \node[vertex] (b7) at (13,1.4) {$b_7$};
        \node[vertex] (a8) at (14.4,0) {$a_8$};
        \node[vertex] (r2) at (-.4+1,-1) {$r_2$};
        \node[vertex] (rp) at (7,-1.3) {$r_{p-1}$};

        \node[vertex] (R') at (11,2.3) {$R'$};
        \node[vertex] (R') at (4,-1.4) {$R$};
        \end{scope}
        
        \tikzstyle{vertex}=[circle, fill=black]
        
        \node[vertex] (a1) at (0,0) {$ $};
        \node[vertex] (b1) at (1,1) {$ $};
        \node[vertex] (a2) at (2,0) {$ $};
        \node[vertex] (b2) at (3,1) {$ $};
        \node[vertex] (a3) at (4,0) {$ $};
        \node[vertex] (b3) at (5,1) {$ $};
        \node[vertex] (a4) at (6,0) {$ $};
        \node[vertex] (b4) at (7,1) {$ $};
        \node[vertex] (a5) at (8,0) {$ $};
        \node[vertex] (b5) at (9,1) {$ $};
        \node[vertex] (a6) at (10,0) {$ $};
        \node[vertex] (b6) at (11,1) {$ $};
        \node[vertex] (a7) at (12,0) {$ $};
        \node[vertex] (b7) at (13,1) {$ $};
        \node[vertex] (a8) at (14,0) {$ $};

        \node[vertex] (r1) at (1,-1) {$ $};
        \node[vertex] (r3) at (7,-1) {$ $};

        \tikzset{decoration={snake,amplitude=.3mm,segment length=3mm, post length=0mm,pre length=0mm}}
        
        \draw(a1)--(b1)--(a2)--(b2)--(a3)--(b3)--(a4)--(b4)--(a5)--(b5)--(a6)--(b6)--(a7)--(b7)--(a8);

        \draw[ForestGreen](a1)--(a2)--(a3)--(a4)--(a5)--(a6)--(a7)--(a8);

        \draw[ForestGreen](r3) to (b4);
        \draw[ForestGreen](b4) to (a8);

        \draw(a1)--(r1);
        \draw[decorate](r1) to (r3);
        \draw(a5)--(r3);
        \draw[rounded corners=3mm] (a8) -- (13.5,-1) --(r3);
        \draw[rounded corners=3mm] (a8) -- (13.5,2) -- (8.5,2) -- (a5);

        \draw[ForestGreen,decorate,bend right=5](a1) to (r3);
            
        \end{tikzpicture}
        \caption{In this example, we follow the assumptions set at the beginning of Case 1 of Lemma \ref{lem:upperBound2Dist} where $k = 4$. Thus,  $\diam(G_2) = 7$, $\ell+1 = 5$, and $k + 4 = 8$. Let the black edges denote $E(G)$ and let the green edges denote $E(G_2)$. Let $b_4 \not\in N_G[r_{p-1}] \cup N_G[a_8]$. The green walk $D_2(R \mathord{-} b \mathord{-} R')$ has a length less than $d_{G_2}(a_1,a_8)$, a contradiction within the proof.}
        \label{fig:Case1_XandB}

\end{figure}

We claim that $r_{p-1} \in N_G[B]$. Suppose not. For a fixed $b \in B$ we have $r_{p-1} \not\in N_G[b]$. Hence, $a_{k+4} \in N[b]$. If $a_{k+4} = b$, then $r_{p-1} \in N_G[b]$ (recall that $r_{p-1} \in N_G(a_{k+4})$), a contradiction. Therefore, $a_{k+4} \neq b$. Now $G$ has a walk $Y = \langle a_{k+4} b a_{\ell} \rangle$ with all vertices distinct. If $Y$ is induced then $D_2(Y)$ gives $d_{G_2}(a_{k+4},a_{\ell}) = 1$, a contradiction. Therefore, $Y$ is not induced, meaning that $a_{k+4} \in N_G(a_{\ell})$ and so $R' \mathord{-} a_{\ell}$ is a path in $G$. The existence of $R' \mathord{-} a_{\ell}$ in $G$ implies that $a_{k+4} \in C_{N_G}(a_{\ell},a_{\ell+4})$ (i.e. $a_{k+4} \in B$). This immediately shows that $r_{p-1} \in N_G[B]$ because $r_{p-1}$ neighbors $a_{k+4} \in R'$, a contradiction. We have proven our claim. We consider two specific subcases of the result of this claim.

\textbf{Case 1.1:} $r_{p-1} \in N_G(B)$ and $r_{p-1} \not\in B$. In particular, let $b' \in B$ and $r_{p-1} \in N_G(b')$. Recall that $r_p = a_{\ell+1}$. If $r_{p-1} = a_{\ell}$, then trivially $a_{\ell} \not\in N_{G_2}(a_{\ell+1})$, a contradiction to the definition of $P_2$. Thus, $r_{p-1} \neq a_{\ell}$.

Next, we will prove that $r_{p-1} \in N_G(a_{\ell})$. Suppose otherwise. Recalling the definition of $k'$, note that $D_2(R)$ has length at most $\lfloor \frac{k'}{2} \rfloor$. There exists an $a_1,a_{\ell}$-walk $Q = D_2(R \mathord{-} b' \mathord{-} a_{\ell})$ of length at most $\lfloor \frac{k'}{2} \rfloor + 1$ in $G_2$. It is easy to verify that $\lfloor \frac{k'}{2} \rfloor + 1 < (\ell - 1) = (\lceil \frac{k+4}{2} \rceil - 1)$ for any $k \geq 3$. Consequently, the length of $Q$ contradicts the assumption that $d_{G_2}(a_1,a_{\ell}) = \ell - 1$. Therefore, $r_{p-1} \in N_G(a_{\ell})$. Now $r_{p-1} \in C_{N_G}(a_{\ell},a_{k+4})$ (i.e. $r_{p-1} \in B$), a contradiction.

\textbf{Case 1.2:} $r_{p-1} \in B$. Note that $r_{p-1} \in N_G(a_{\ell})$ by the definition of $B$. We also have that $a_{\ell} \in N_G(a_{k+4})$, otherwise $D_2(\langle a_{\ell} r_{p-1} a_{k+4} \rangle)$ gives $d_{G_2}(a_{\ell},a_{k+4}) = 1$, a contradiction. Recall that $R$ has length $3$ or greater. Thus, $r_{p-2}$ and $r_{p-3}$ exist. If $r_{p-3} \in N_G(a_{k+4})$, then $\langle r_{p-3} a_{k+4} a_{\ell+1} \rangle$ has length at most $2$ in $G$, a contradiction to the fact that $R$ is a shortest path with $d_G(r_{p-3},r_p) = 3$. Hence, $r_{p-3} \not\in N_G(a_{k+4})$. See Fig. \ref{fig:Case111-112} for a graph depicting the current state of this proof where $k = 4$.

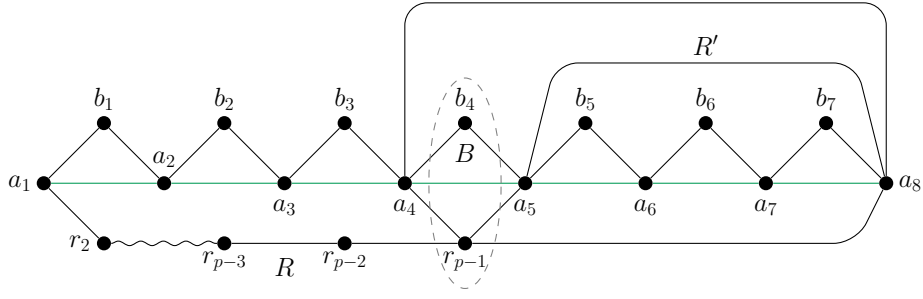
\begin{figure}[h]
    \centering
        \begin{tikzpicture}[scale=.8, every node/.style={scale=.48}]

        \tikzstyle{vertex}=[circle, fill=white]
        \begin{scope}\LARGE
        \node[vertex] (a1) at (-.4,0) {$a_1$};
        \node[vertex] (b1) at (1,1.4) {$b_1$};
        \node[vertex] (a2) at (2,.4) {$a_2$};
        \node[vertex] (b2) at (3,1.4) {$b_2$};
        \node[vertex] (a3) at (4,-.4) {$a_3$};
        \node[vertex] (b3) at (5,1.4) {$b_3$};
        \node[vertex] (a4) at (6,-.4) {$a_4$};
        \node[vertex] (b4) at (7,1.4) {$b_4$};
        \node[vertex] (a5) at (8,-.4) {$a_5$};
        \node[vertex] (b5) at (9,1.4) {$b_5$};
        \node[vertex] (a6) at (10,-.4) {$a_6$};
        \node[vertex] (b6) at (11,1.4) {$b_6$};
        \node[vertex] (a7) at (12,-.4) {$a_7$};
        \node[vertex] (b7) at (13,1.4) {$b_7$};
        \node[vertex] (a8) at (14.4,0) {$a_8$};
        \node[vertex] (r2) at (-.4+1,-1) {$r_2$};
        \node[vertex] (rp) at (7,-1.3) {$r_{p-1}$};
        \node[vertex] (rp) at (5,-1.3) {$r_{p-2}$};
        \node[vertex] (rp) at (3,-1.3) {$r_{p-3}$};
        \node[vertex] (b4) at (7,.5) {$B$};

        \node[vertex] (R') at (11,2.3) {$R'$};
        \node[vertex] (R') at (4,-1.4) {$R$};
        \end{scope}
        
        \tikzstyle{vertex}=[circle, fill=black]
        
        \node[vertex] (a1) at (0,0) {$ $};
        \node[vertex] (b1) at (1,1) {$ $};
        \node[vertex] (a2) at (2,0) {$ $};
        \node[vertex] (b2) at (3,1) {$ $};
        \node[vertex] (a3) at (4,0) {$ $};
        \node[vertex] (b3) at (5,1) {$ $};
        \node[vertex] (a4) at (6,0) {$ $};
        \node[vertex] (b4) at (7,1) {$ $};
        \node[vertex] (a5) at (8,0) {$ $};
        \node[vertex] (b5) at (9,1) {$ $};
        \node[vertex] (a6) at (10,0) {$ $};
        \node[vertex] (b6) at (11,1) {$ $};
        \node[vertex] (a7) at (12,0) {$ $};
        \node[vertex] (b7) at (13,1) {$ $};
        \node[vertex] (a8) at (14,0) {$ $};

        \node[vertex] (r1) at (1,-1) {$ $};
        \node[vertex] (r15) at (3,-1) {$ $};
        \node[vertex] (r2) at (5,-1) {$ $};
        \node[vertex] (r3) at (7,-1) {$ $};

        \tikzset{decoration={snake,amplitude=.3mm,segment length=3mm, post length=0mm,pre length=0mm}}
        
        \draw(a1)--(b1)--(a2)--(b2)--(a3)--(b3)--(a4)--(b4)--(a5)--(b5)--(a6)--(b6)--(a7)--(b7)--(a8);

        \draw[ForestGreen](a1)--(a2)--(a3)--(a4)--(a5)--(a6)--(a7)--(a8);

        \draw(r3) to (a4);
        \draw(a1)--(r1);
        \draw[decorate](r1) to (r15);
        \draw(r15)--(r2);
        \draw(r2)--(r3);
        \draw(a5)--(r3);
        
        \draw[rounded corners=3mm] (a8) -- (13.5,2) -- (8.5,2) -- (a5);
        \draw[rounded corners=3mm] (a8) -- (14,3) -- (6,3) -- (a4);
        \draw[rounded corners=3mm] (a8) -- (13.5,-1) --(r3);

        \draw[dashed,gray] (7,0) ellipse (.6cm and 1.75cm); 
            
        \end{tikzpicture}
        \caption{This graph shows the setup of Lemma \ref{lem:upperBound2Dist} after entering Case 1.2 where $k = 4$. The dotted ellipse represents $B$. Black edges belong to $E(G)$ and green edges belong to $E(G_2)$.}
        \label{fig:Case111-112}

\end{figure}

If $r_{p-2} \not\in N_G(a_{k+4}) \cup N_G(a_{\ell})$, then there exists $Q = \langle a_{k+4} r_{p-1} r_{p-2} r_{p-1} a_{\ell} \rangle$ in $G$; in fact, $Q$ exists even if $r_{p-2} \in \{a_{\ell},a_{k+4}\}$. Now $D_2(Q) = \langle a_{k+4} r_{p-2} a_{\ell} \rangle$ is an $a_{k+4},a_{\ell}$-walk of length at most $2$ in $G_2$, a contradiction. We have shown that $r_{p-2} \in N_G(a_{k+4}) \cup N_G(a_{\ell})$. In particular, we will show that $r_{p-2} \in N_G(a_{\ell})$. Suppose not, and thus $r_{p-2} \in N_G(a_{k+4})$. Now $D_2(R[r_1,r_{p-2}] \mathord{-} a_{k+4})$ is an $a_1,a_{k+4}$-path of length at most $\lfloor \frac{k}{2} \rfloor$ in $G_2$, a contradiction. Therefore, $r_{p-2} \in N_G(a_{\ell})$.

There exists $Q' = R[r_1,r_{p-2}] \mathord{-} a_{\ell}$ in $G$. If $r_{p-3} \not\in N_G(a_{\ell})$, then $D_2(Q')$ is an $a_1,a_{\ell}$-path of length at most $\lfloor \frac{k'}{2} \rfloor$. Then, noting that $\lfloor \frac{k'}{2} \rfloor < (\ell - 1) = (\lceil \frac{k+4}{2} \rceil - 1)$ for any $k \geq 3$, observe that $D_2(Q')$ is shorter than $d_{G_2}(a_1,a_{\ell}) = \ell - 1$, a contradiction. Thus, $r_{p-3} \in N_G(a_{\ell})$.

There exists $T = R[r_1,r_{p-3}] \mathord{-} a_{\ell} \mathord{-} a_{k+4}$ in $G$. Recalling that $r_{p-3} \not\in N_G(a_{k+4})$, observe that $D_2(T)$ is an $a_1,a_{k+4}$-path of length at most $\lfloor \frac{k}{2} \rfloor$ in $G_2$, a contradiction.

\textbf{Case 2:} $r_{p-1}$ is \emph{not} in the closed neighborhood of every vertex in $R'$ in $G$. Let $r'_f$ be the first vertex in $R'$ that is not in $N_G[r_{p-1}]$. Since $r_{p-1}$ is adjacent to $a_{\ell+1} = r'_1$ by the definitions of $R$ and $R'$, clearly $f > 1$. There exists $T = R[r_1,r_{p-1}] \mathord{-} R'[r'_{f-1},r'_{p'}]$ in $G$. Notice that $T$ is an $a_1,a_{k+4}$-walk with length no greater than $2k$.

It is straightforward to check that if $f$ is even, then $T$ has even length. Subsequently, from $D_2(T)$ we can extract an $a_1,a_{k+4}$-path having length at most $k$. This contradicts the assumption that $d_{G_2}(a_1,a_{k+4}) = k + 3$. So, $f$ is odd. See Fig. \ref{fig:Case2}.

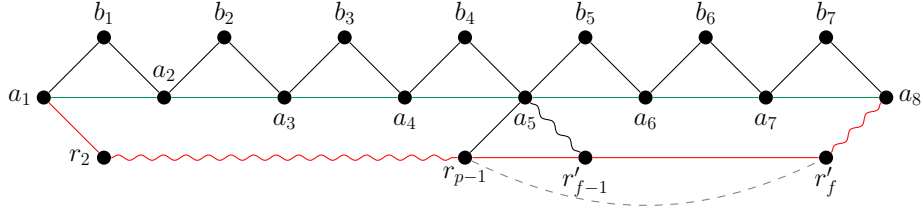
\begin{figure}[h]
    \centering
        \begin{tikzpicture}[scale=.8, every node/.style={scale=.48}]

        \tikzstyle{vertex}=[circle, fill=white]
        \begin{scope}\LARGE
        \node[vertex] (a1) at (-.4,0) {$a_1$};
        \node[vertex] (b1) at (1,1.4) {$b_1$};
        \node[vertex] (a2) at (2,.4) {$a_2$};
        \node[vertex] (b2) at (3,1.4) {$b_2$};
        \node[vertex] (a3) at (4,-.4) {$a_3$};
        \node[vertex] (b3) at (5,1.4) {$b_3$};
        \node[vertex] (a4) at (6,-.4) {$a_4$};
        \node[vertex] (b4) at (7,1.4) {$b_4$};
        \node[vertex] (a5) at (8,-.4) {$a_5$};
        \node[vertex] (b5) at (9,1.4) {$b_5$};
        \node[vertex] (a6) at (10,-.4) {$a_6$};
        \node[vertex] (b6) at (11,1.4) {$b_6$};
        \node[vertex] (a7) at (12,-.4) {$a_7$};
        \node[vertex] (b7) at (13,1.4) {$b_7$};
        \node[vertex] (a8) at (14.4,0) {$a_8$};
        \node[vertex] (r2) at (-.4+1,-1) {$r_2$};
        \node[vertex] (rp) at (7,-1.3) {$r_{p-1}$};

        \node[vertex] (sr'f-1) at (9,-1.4) {$r'_{f-1}$};
        \node[vertex] (sr'f) at (13,-1.4) {$r'_f$};
        \end{scope}
        
        \tikzstyle{vertex}=[circle, fill=black]
        
        \node[vertex] (a1) at (0,0) {$ $};
        \node[vertex] (b1) at (1,1) {$ $};
        \node[vertex] (a2) at (2,0) {$ $};
        \node[vertex] (b2) at (3,1) {$ $};
        \node[vertex] (a3) at (4,0) {$ $};
        \node[vertex] (b3) at (5,1) {$ $};
        \node[vertex] (a4) at (6,0) {$ $};
        \node[vertex] (b4) at (7,1) {$ $};
        \node[vertex] (a5) at (8,0) {$ $};
        \node[vertex] (b5) at (9,1) {$ $};
        \node[vertex] (a6) at (10,0) {$ $};
        \node[vertex] (b6) at (11,1) {$ $};
        \node[vertex] (a7) at (12,0) {$ $};
        \node[vertex] (b7) at (13,1) {$ $};
        \node[vertex] (a8) at (14,0) {$ $};

        \node[vertex] (r1) at (1,-1) {$ $};
        \node[vertex] (r3) at (7,-1) {$ $};

        \node[vertex] (r'f-1) at (9,-1) {$ $};
        \node[vertex] (r'f) at (13,-1) {$ $};

        \tikzset{decoration={snake,amplitude=.35mm,segment length=3mm, post length=0mm,pre length=0mm}}
        
        \draw(a1)--(b1)--(a2)--(b2)--(a3)--(b3)--(a4)--(b4)--(a5)--(b5)--(a6)--(b6)--(a7)--(b7)--(a8);

        \draw[ForestGreen](a1)--(a2)--(a3)--(a4)--(a5)--(a6)--(a7)--(a8);

        \draw[red](a1)--(r1);
        \draw[decorate,red](r1) to (r3);
        \draw(a5)--(r3);

        \draw[decorate](a5) to (r'f-1);
        \draw[red](r'f-1)--(r'f);
        \draw[decorate,red](r'f) to (a8);

        \draw[red](r3) to (r'f-1);
        \draw[bend right=26,dashed,gray](r3) to (r'f);
            
        \end{tikzpicture}
        \caption{This graph shows the setup of Lemma \ref{lem:upperBound2Dist} immediately prior to entering Cases 2.1 or 2.2. The black and red edges belong to $E(G)$ while green edges belong to $E(G_2)$. The dotted arc is an edge that does not exist in $G$. The red path denotes $T$.}
        \label{fig:Case2}

\end{figure}

\textbf{Case 2.1:} $f > 3$. Hence, $f \geq 5$ (since $f$ cannot be even) and $R'$ has length $5$ or greater. There exists a walk $R'[r'_{p'},r'_{f-1}] \mathord{-} \langle r_{p-1} r_{p} \rangle$ in $G$ whose length is strictly less than the length of $R'$, a contradiction to the assumption that $R'$ is a shortest path in $G$.

\textbf{Case 2.2:} $f = 3$. So, $r'_3 = r'_f$ and $r'_2 = r'_{f-1}$. Notice that $\langle r_{p-1} r'_2 r'_3 \rangle$ is induced in $G$. There exists $Q = R[r_1,r_{p-1}] \mathord{-} \langle r'_2 r'_3 r'_2 r'_1\rangle$ in $G$. It follows that $D_2(Q)$ is an $a_1,a_{\ell+1}$-walk with length at most $\tau = \lfloor \frac{k'}{2} \rfloor + 2$ in $G_2$. It is easy to verify that $\tau < \ell$ for any value of $k \geq 3$, a contradiction to the assumption that $d_{G_2}(a_1,a_{\ell+1}) = \ell$.
\end{proof}

The following simple result will be useful in later proofs.

\begin{proposition}
Given a graph $G$ with $\diam(G) = k \geq 3$, let $P_2 = \langle a_1 a_2 \cdots a_{k+3} a_{k+4} \rangle$ be a shortest path in $G_2$, and let $i \in [3,k+2]$ such that $a_1,a_{k+4} \in N_G(a_i)$. The set $\{a_1,a_{i},a_{k+4}\}$ induces a triangle in $G$.
\label{prop:centerTriangle}
\end{proposition}

\begin{proof}
There exists a path $R = \langle a_1 a_{i} a_{k+4} \rangle$ in $G$. If $a_1 \not\in N_G(a_{k+4})$, then $D_2(R) = \langle a_1 a_{k+4} \rangle$ exists in $G_2$, a contradiction to the assumption that $P_2$ is a shortest path.
\end{proof}

Next, we will show that any central vertex in $P_2$ that is adjacent to both the first and last vertices in $P_2$ infers the existence of two additional edges in $G$.

\begin{lemma}
Under the hypothesis of Proposition \ref{prop:centerTriangle} where $i > 3$, we have that $a_1,a_{k+4} \in N_G(a_{i-1})$.
\label{lem:nextForcedEdges}
\end{lemma}

\begin{proof}
Let $b \in C_{N_G}(a_{i-1},a_i)$ and $R = \langle a_1 a_i b a_i a_{k+4} \rangle$. By Proposition \ref{prop:centerTriangle} we have that $a_1 \in N_G(a_{k+4})$. Suppose for the sake of contradiction that $b \not\in N_G(a_1) \cup N_G(a_{k+4})$. Then, there exists a walk $D_2(R) = \langle a_1 b a_{k+4} \rangle$ (this walk exists even if $b \in \{a_1,a_{k+4}\}$). But, now we have that $d_{G_2}(a_1,a_{k+4}) \leq 2$, a contradiction. Thus, $b \in N_G(a_1) \cup N_G(a_{k+4})$.

Next, if $b \in N_G(a_1)$, then let $(u,v) = (a_1,a_{k+4})$; alternatively, if $b \in N_G(a_{k+4})$ then let $(u,v) = (a_{k+4},a_1)$. Surely $u \in N_G(a_{i-1})$ when $b = u$ (trivially) or when $b \neq u$ (otherwise $D_2(\langle u b a_{i-1} \rangle) = \langle u a_{i-1} \rangle$ exists in $G_2$, a contradiction). Furthermore, we have $v \in N_G(a_{i-1})$ when $b = v$ (trivially) or when $b \neq v$ (otherwise $D_2(\langle v u a_{i-1} \rangle) = \langle v a_{i-1} \rangle$ exists in $G_2$, a contradiction). We have shown that $u,v \in N_G(a_{i-1})$, so we are done.
\end{proof}

Next, we will show that a graph satisfying the conclusions of Lemmas \ref{lem:upperBound2Dist} and \ref{lem:nextForcedEdges} contains a subgraph that induces the complement of a path. The existence of this structure will be vital for proving Theorem \ref{lem:k+2Upper}.

\begin{theorem}
Let $G$ have $\diam(G) = k \geq 3$ and let $P_2 = \langle a_1 a_2 \cdots a_{k+3} a_{k+4} \rangle$ be a shortest path in $G_2$. Then, $V(P_2)$ induces the complement of a path in $G$.
\label{lem:AInducedComp}
\end{theorem}

\begin{proof}
For easier reference, let $A = V(P_2)$. Let $\ell = \lceil \frac{k+4}{2} \rceil$ and let $P_2^R = \langle a_{k+4} a_{k+3} \cdots a_1 \rangle$ be a path in $G_2$ (that is, $P_2^R$ is the reverse of $P_2$). By Lemma \ref{lem:upperBound2Dist} we have that $a_1,a_{k+4} \in N_G(a_{\ell+1})$. By Proposition \ref{prop:centerTriangle} where $i := \ell$ we have that $a_1 \in N_G(a_{k+4})$.

We will show by induction that, for each $3 \leq j \leq \ell$, we have $a_1,a_{k+4} \in N_G(a_{j})$. The statement is true if $j = \ell$ because Lemma \ref{lem:nextForcedEdges} is applicable by setting $P_2 := P_2$ and $i := \ell+1$, giving $a_1,a_{k+4} \in N_G(a_{\ell})$. Now suppose that $3 \leq j < \ell$. By the inductive hypothesis, $a_1,a_{k+4} \in N_G(a_{j+1})$. Then, Lemma \ref{lem:nextForcedEdges} holds where $P_2 := P_2$, and $i := j+1$, promising that $a_1,a_{k+4} \in N_G(a_{j})$. Fig. \ref{fig:goingLeft} demonstrates the conclusion of our inductive statement.

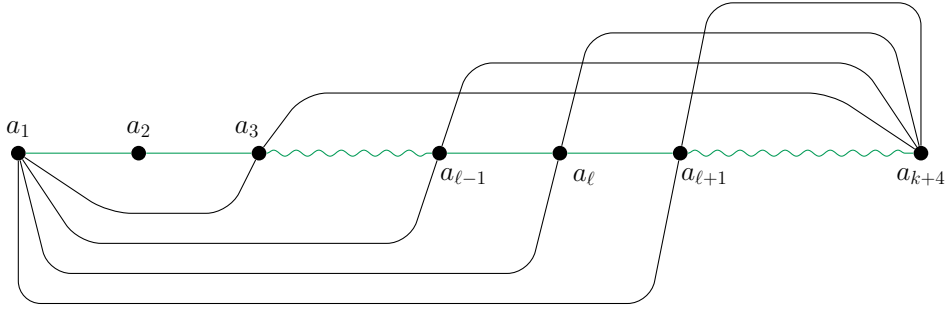
\begin{figure}[h]
    \centering
        \begin{tikzpicture}[scale=.8, every node/.style={scale=.5}]

        \tikzstyle{vertex}=[circle, fill=white]
        \begin{scope}\LARGE
        \node[vertex] (a1) at (0,0.4) {$a_1$};
        \node[vertex] (a2) at (2,.4) {$a_2$};
        \node[vertex] (a3) at (4-.2,.4) {$a_3$};
        \node[vertex] (a) at (7.4,-.4) {$a_{\ell-1}$};
        \node[vertex] (a5) at (9.4,-.4) {$a_{\ell}$};
        \node[vertex] (a6) at (11.4,-.4) {$a_{\ell+1}$};
        \node[vertex] (a8) at (15,-.4) {$a_{k+4}$};
        \end{scope}
        
        \tikzstyle{vertex}=[circle, fill=black]
        
        \node[vertex] (a1) at (0,0) {$ $};
        \node[vertex] (a2) at (2,0) {$ $};
        \node[vertex] (a3) at (4,0) {$ $};
        \node[vertex] (a4) at (7,0) {$ $};
        \node[vertex] (a5) at (9,0) {$ $};
        \node[vertex] (a6) at (11,0) {$ $};
        \node[vertex] (a8) at (15,0) {$ $};

        \tikzset{decoration={snake,amplitude=.4mm,segment length=3mm, post length=0mm,pre length=0mm}}
        
        \draw[ForestGreen](a1)--(a2)--(a3);
        \draw[ForestGreen,decorate](a3)--(a4);
        \draw[ForestGreen](a4)--(a5);
        \draw[ForestGreen](a5)--(a6);
        \draw[ForestGreen,decorate](a6)--(a8);

        \draw[rounded corners=3mm] (a8) -- (15,2.5) -- (11.5,2.5) -- (a6);
        \draw[rounded corners=3mm] (a1) -- (0,-2.5) -- (10.5,-2.5) -- (a6);

        \draw[rounded corners=3mm] (a8) -- (14.5,2) -- (9.5,2) -- (a5);
        \draw[rounded corners=3mm] (a1) -- (.5,-2) -- (8.5,-2) -- (a5);
        
        \draw[rounded corners=3mm] (a4) -- (7.5,1.5) -- (14,1.5) -- (a8);
        \draw[rounded corners=3mm] (a1) -- (1,-1.5) -- (6.5,-1.5) -- (a4);

        \draw[rounded corners=3mm] (a3) -- (4.75,1) -- (13.5,1) -- (a8);
        \draw[rounded corners=3mm] (a1) -- (1.5,-1) -- (3.5,-1) -- (a3);
        \end{tikzpicture}
        \caption{This graph, where black edges and green edges belong to $E(G)$ and $E(G_2)$, respectively, depicts an example of the edges provided via repeated application of Lemma \ref{lem:nextForcedEdges}. As a result, for each $3 \leq j \leq \ell+1$ we have $a_1,a_{k+4} \in N_G(a_j)$.}
        \label{fig:goingLeft}

\end{figure}

Next, we will also show by induction that, for each $\ell \leq j \leq k+2$, we have $a_1,a_{k+4} \in N_G(a_{j})$. The statement is true if $j = \ell$ because Lemma \ref{lem:nextForcedEdges} is applicable by setting $P_2 := P_2^R$ and $i := \ell-1$, giving $a_1,a_{k+4} \in N_G(a_{\ell})$. Now suppose that $\ell < j \leq k+2$. By the inductive hypothesis, $a_1,a_{k+4} \in N_G(a_{j-1})$. Then, Lemma \ref{lem:nextForcedEdges} holds where $P_2 := P_2^R$ and $i := j-1$, promising that $a_1,a_{k+4} \in N_G(a_{j})$. We have proven our claim.

Observe now that both $a_1$ and $a_{k+4}$ are in the closed neighborhood of every vertex in $A \setminus \{a_2,a_{k+3}\}$. We claim that for each $i,j \in [3,k + 2]$ such that $\abs(i - j) > 1$ we have $a_i \in N_G(a_j)$. Suppose on the contrary that for some $i',j' \in [3,k+2]$ where $\abs(i' - j') > 1$ we have $a_{i'} \not\in N_G(a_{j'})$. The existence of $D_2(\langle a_{i'} a_1 a_{j'} \rangle) = \langle a_{i'} a_{j'} \rangle$ gives $d_{G_2}(a_{i'},a_{j'}) = 1$, a contradiction. We have proven the claim. At this stage of the proof, $A \setminus \{a_2, a_{k+3}\}$ has been shown to induce the complement of a path in $G$.

We claim that $a_2 \in N_G(a_{k+4})$ and $a_{k+3} \in N_G(a_1)$. Note the symmetry of the fact that $a_2$ and $a_{k+3}$ are the second and second to last vertices in $P_2$, respectively. W.l.o.g. suppose for the sake of contradiction that $a_2 \not\in N_G(a_{k+4})$. Note that if $d_G(a_2,a_{k+4}) = 2$, then $d_{G_2}(a_2,a_{k+4}) = 1$, a contradiction. There exists a vertex $b_1 \in C_{N_G}(a_1,a_2)$ because we have assumed that $d_{G_2}(a_1,a_2) = 1$. Let $W = \langle a_{k+4} a_1 b_1 a_1 a_{\ell} \rangle$ in $G$. If $b_1 \not\in N_G(a_{\ell}) \cup N_G(a_{k+4})$, then $D_2(W) = \langle a_{k+4} b_1 a_{\ell} \rangle$ is shorter than the assumed value of $d_{G_2}(a_{\ell},a_{k+4})$, a contradiction (this is easily verified for any $k$, and is true even if $b_1 \in \{a_{\ell},a_{k+4}\}$). But, $b_1 \in N_G(a_{k+4})$ gives $d_G(a_2,a_{k+4}) = 2$, a contradiction. Hence, it is necessary that $b_1 \in N_G(a_{\ell})$. However, the existence of $\langle a_2 b_1 a_{\ell} \rangle$ in $G$ implies that $d_{G_2}(a_2, a_{\ell}) = 1$, a contradiction. Thus, $a_2 \in N_G(a_{\ell})$. Because $a_2 \not\in N_G(a_{k+4})$, we have that the existence of $\langle a_2 a_{\ell} a_{k+4} \rangle$ in $G$ gives $d_{G_2}(a_2,a_{k+4}) = 1$, a contradiction. We have proven our claim.

Note the symmetry of the fact that $a_2$ and $a_{k+3}$ are adjacent to $a_{k+4}$ and $a_1$, respectively. We will prove that $a_2 \in N_G(u)$ for all $u \in A \setminus \{a_1,a_2,a_3,a_{k+3},a_{k+4}\}$. This will similarly show that $a_{k+3} \in N_G(v)$ for all $v \in A \setminus \{a_{k+4},a_{k+3},a_{k+2},a_2,a_1\}$. W.l.o.g. suppose for the sake of contradiction that $a_2 \not\in N_G(u')$ for some $u' \in A \setminus \{a_1,a_2,a_3,a_{k+3},a_{k+4}\}$. It follows that $\langle a_2 a_{k+4} u' \rangle$ is an induced path in $G$ that gives $d_{G_2}(a_2,u') = 1$, a contradiction. We have proven our claim.

In order to complete the proof, it remains to show that $a_2 \in N_G(a_{k+3})$. If not, then the induced path $\langle a_2 a_4 a_{k+3} \rangle$ in $G$ gives $d_{G_2}(a_2,a_{k+3}) = 1$, a contradiction.
\end{proof}

The next lemma is crucial because it shows that the existence of $V(P_2)$, the complement of a path in $G$, applies special constraints on its neighboring vertices. More specifically, any such vertex has at least two edges between itself and $V(P_2)$, and these edges are incident to vertices that are far apart in $P_2$.
      
\begin{lemma}
Given a graph $G$ with $\diam(G) = k \geq 3$, let $P_2 = \langle a_1 a_2 \cdots a_{k+3} a_{k+4} \rangle$ be a shortest path in $G_2$ and let $A \subseteq G$ be the subgraph induced by $V(P_2)$. Let $u \in N_G(A)$. Then, $N_G(u) \cap A = \{a_i,a_j\}$ such that $i,j \in [1,k+4]$ and $\abs(i - j) > 2$.
\label{lem:compPathLowDiam}
\end{lemma}

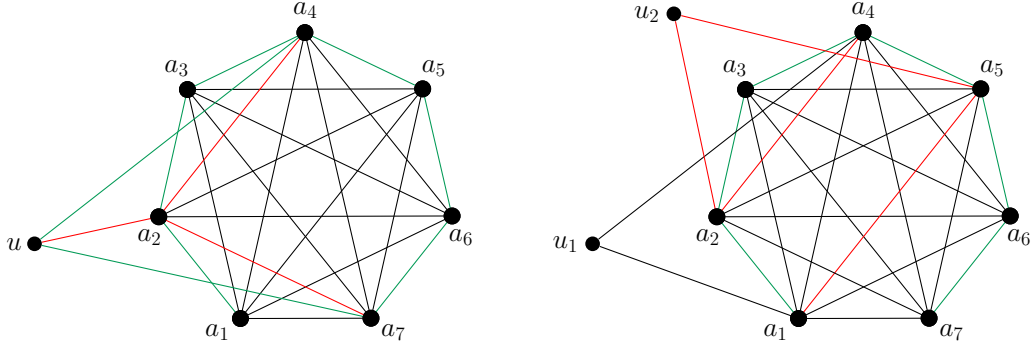
\begin{figure}[h]
    \centering

        \begin{tikzpicture}[scale=.9, every node/.style={scale=.5}]

        \tikzstyle{vertex}=[circle, fill=white]
        \begin{scope}\LARGE
            \node[vertex] (u) at (-4-.3,-.9) {$u$};
            
            \node[vertex] (a1) at (-1.3,-2.2) {$a_1$};
            \node[vertex] (a7) at (1.3,-2.2) {$a_7$};
            
            \node[vertex] (a2) at (-2.3,-.8) {$a_2$};
            \node[vertex] (a6) at (2.3,-.8) {$a_6$};

            \node[vertex] (a3) at (-1.9,1.7) {$a_3$};
            \node[vertex] (a5) at (1.9,1.7) {$a_5$};
            
            \node[vertex] (a4) at (0,2.55) {$a_4$};
        \end{scope}

        \tikzstyle{vertex}=[circle, fill=black]
        \tkzDefPoint(0,0){O}\tkzDefPoint(.97,-2){A}
        \tkzDefPointsBy[rotation=center O angle 360/7](A,B,C,D,E,F){B,C,D,E,F,G}
        \tkzDrawPoints[fill =black,size=12,color=black](A,B,C,D,E,F,G)
        \draw[ForestGreen](A)--(B)--(C)--(D)--(E)--(F)--(G);
        
        \node[vertex] (v) at (-4,-.9) {$ $};
        
        \tikzset{decoration={snake,amplitude=.3mm,segment length=3mm, post length=0mm,pre length=0mm}}
        \draw(A)edge(C)edge(D)edge(E)edge(G);
        \draw(B)edge(D)edge(E)edge(F)edge(G);
        \draw(C)edge(E)edge(F)edge(G);
        \draw(D)edge(G);
        \draw(E)edge(G);

        \draw[red](A)--(F);
        \draw[red](D)--(F);
        \draw[red](v)--(F);
        \draw[ForestGreen](v)edge(D);
        \draw[ForestGreen](v)edge(A);

        \tkzDefPoint(0,0){O}\tkzDefPoint(.97,-2){A}
        \tkzDefPointsBy[rotation=center O angle 360/7](A,B,C,D,E,F){B,C,D,E,F,G}
        \tkzDrawPoints[fill =black,size=12,color=black](A,B,C,D,E,F,G)
 
        \end{tikzpicture}
        ~~~
        \begin{tikzpicture}[scale=.9, every node/.style={scale=.5}]

        \tikzstyle{vertex}=[circle, fill=white]
        \begin{scope}\LARGE
            \node[vertex] (u) at (-4-.4,-.9) {$u_1$};
            
            \node[vertex] (a1) at (-1.3,-2.2) {$a_1$};
            \node[vertex] (a7) at (1.3,-2.2) {$a_7$};
            
            \node[vertex] (a2) at (-2.3,-.8) {$a_2$};
            \node[vertex] (a6) at (2.3,-.8) {$a_6$};

            \node[vertex] (a3) at (-1.9,1.7) {$a_3$};
            \node[vertex] (a5) at (1.9,1.7) {$a_5$};
            
            \node[vertex] (a4) at (0,2.55) {$a_4$};
            \node[vertex] (v2') at (-2.8-.4,2.5) {$u_2$};
        \end{scope}

        \tikzstyle{vertex}=[circle, fill=black]
        \tkzDefPoint(0,0){O}\tkzDefPoint(.97,-2){A}
        \tkzDefPointsBy[rotation=center O angle 360/7](A,B,C,D,E,F){B,C,D,E,F,G}
        \tkzDrawPoints[fill =black,size=12,color=black](A,B,C,D,E,F,G)
        \draw[ForestGreen](A)--(B)--(C)--(D)--(E)--(F)--(G);
        
        \node[vertex] (v) at (-4,-.9) {$ $};
        \node[vertex] (v2) at (-2.8,2.5) {$ $};
        
        \tikzset{decoration={snake,amplitude=.3mm,segment length=3mm, post length=0mm,pre length=0mm}}
        \draw(A)edge(C)edge(D)edge(E)edge(G);
        \draw(B)edge(D)edge(E)edge(F)edge(G);
        \draw(C)edge(E);
        \draw[red](C)edge(G);
        \draw(D)edge(G);
        \draw(E)edge(G);
        \draw(C)edge(F);

        \draw(A)--(F);
        \draw[red](D)--(F);
        \draw(v)edge(D);
        \draw(v)edge(G);
        \draw[red](v2)edge(C);
        \draw[red](v2)edge(F);

        \tkzDefPoint(0,0){O}\tkzDefPoint(.97,-2){A}
        \tkzDefPointsBy[rotation=center O angle 360/7](A,B,C,D,E,F){B,C,D,E,F,G}
        \tkzDrawPoints[fill=black,size=12,color=black](A,B,C,D,E,F,G)
        \end{tikzpicture}
        
        \caption{The black and red edges belong to $E(G)$ and the green edges belong to $E(G_2)$. On the left, the graph shows a contradiction that occurs within the proof of Lemma \ref{lem:compPathLowDiam} where $(k,h,s,t) = (3,2,4,7)$. The red edges, in particular, denote $W$. On the right, we portray the proof of Theorem \ref{lem:k+2Upper} where $k = 3$, $(b,d) = (1,4)$, and $(p,q) = (2,5)$. The red edges denote $W$ within Case $1$. Note that the path $D_2(W) = \langle a_1 u_2 a_4 \rangle$ exists but is not drawn.}
        \label{fig:ComplementPath+v}

\end{figure}

\begin{proof}
We will prove the statement directly. By Theorem \ref{lem:AInducedComp}, $A$ is the complement of a path. Let $u \in N_G(a_h)$ for any $h$. If $h = 1$, then let $S = \{a_h,a_{h+1}\}$; otherwise if $h = k+4$, then let $S = \{a_{h-1},a_h\}$; otherwise, let $S = \{a_{h-1},a_h,a_{h+1}\}$. Next, let $a_s,a_t \in A \setminus S$ such that $i \not\in \{s,t\}$ and $\abs(s-t) > 2$ (it is easy to verify that such a pair $a_s,a_t$ exists). There exists a walk $W = \langle a_s a_h u a_h a_t \rangle$ in $G$. If $u \not\in N_G(a_s) \cup N_G(a_t)$, then $D_2(W) = \langle a_s u a_t \rangle$ exists in $G_2$, giving $d_{G_2}(a_s,a_t) = 2$, a contradiction to the assumption that $\abs(s - t) > 2$. See Fig. \ref{fig:ComplementPath+v} for an example of this contraction. Therefore, $u \in N_G(a_s) \cup N_G(a_t)$. If $u \in N_G(a_s)$ and $\abs(h - s) > 3$, then we are done by setting $(i,j) := (h,s)$. If $u \in N_G(a_t)$ and $\abs(h - t) > 3$, then we are done by setting $(i,j) := (h,t)$. W.l.o.g. it remains for us to find appropriate $(i,j)$ when $u \in N_G(a_t)$ and $\abs(h - t) = 2$.

To simplify, there exists $N_G(u) \cap A = \{a_h,a_t\}$ such that $\abs(h - t) = 2$. W.l.o.g. suppose that $h < t$. Clearly, $h \in [1,k+2]$. We consider two cases based on the value of $h$.

\textbf{Case 1:} $1 < h < k + 2$. Note that $d_{G_2}(a_{h-1},a_{t+1}) = 4$. There exists a path $T = \langle a_{h-1}  a_t u a_h a_{t+1} \rangle$ in $G$. If $u \not\in N_G(a_{h-1}) \cup N_G(a_{t+1})$, then $D_2(T) = \langle a_{h-1} u a_{t+1} \rangle$ gives $d_{G_2}(a_{h-1},a_{t+1}) = 2$, a contradiction. Thus, $u \in N_G(a_{h-1}) \cup N_G(a_{t+1})$. If $u \in N_G(a_{h-1})$, then we are done by setting $(i,j) := (h-1,t)$. But if $u \in N_G(a_{t+1})$, then we are done by setting $(i,j) := (h,t+1)$.

\textbf{Case 2:} $h = 1$ or $h = k + 2$. More specifically, this means that either $(h,t) = (1,3)$ or $(h,t) = (k+2,k+4)$. Due to the symmetry, finding indices that give us appropriate $(i,j)$ in one of the two cases is enough to complete the proof. Hence, w.l.o.g. we may assume that $(h,t) = (1,3)$. We will show that $u \in N_G(a_4) \cup N_G(a_7)$. Suppose not. There exists a path $R = \langle a_4 a_1 u a_3 a_7 \rangle$ in $G$. Then, $D_2(R) = \langle a_4 u a_7 \rangle$ gives $d_{G_2}(a_4,a_7) = 2$, a contradiction. Thus, $u \in N_G(a_4) \cup N_G(a_7)$. In particular, if $u \in N_G(a_4)$, then we are done by setting $(i,j) := (1,4)$. But if $u \in N_G(a_7)$, then we are done by setting $(i,j) := (1,7)$.
\end{proof}

We are prepared to prove our main result:

\begin{theorem}
Given a graph $G = (V,E)$ with $\diam(G) = k \geq 3$, we have that either $G_2$ is disconnected or $\diam(G_2) \leq k + 2$.
\label{lem:k+2Upper}
\end{theorem}

\begin{proof}
Suppose that $G_2$ is connected and $\diam(G_2) > k + 2$. Thus, $G_2$ contains a diametral path of length $k + 3$ or greater. This path necessarily contains a shortest path $P_2 = \langle a_1 a_2 \cdots a_{k+3} a_{k+4} \rangle$. By Theorem \ref{lem:AInducedComp}, we have that $V(P_2)$ induces the complement of a path in $G$. Let $A \subseteq G$ be the subgraph induced by $V(P_2)$ for simplicity. It is trivial to see that $\diam(A) = 2$. Since $\diam(G) = k$, there exists $S = N_G(A)$ where $S \neq \emptyset$.

We claim that $G \setminus (A \cup S) = \emptyset$. Suppose not. Let $u \in S$ and $v \in N_G(u) \setminus (A \cup S)$. By Lemma \ref{lem:compPathLowDiam}, $u \in S$ satisfies $N_G(u) \cap A = \{a_i,a_j\}$ such that $i,j \in [1,k+4]$ and $\abs(i - j) > 2$. Clearly, for any $a \in A$ we have $d_G(u,a) \leq 2$. Consequently, there exists $T = \langle a_i u v u a_j \rangle$ in $G$. If $v \not\in N_G(a_i) \cup N_G(a_j)$, then $D_2(T) = \langle a_i v a_j \rangle$ in $G_2$ gives $d_{G_2}(a_i,a_j) = 2$, a contradiction to the fact that $\abs(i - j) > 2$. Thus, $v \in N_G(a_i) \cup N_G(a_j)$, a contradiction to the fact that $v \not\in S$. We have proven the claim. 

It is easy to see that $|S| > 1$, otherwise $\diam(G[V \cup S]) = 2$, a contradiction. Next, we will prove that for all pairs $v_1,v_2 \in S$ we have $d_G(v_1,v_2) \leq 2$. On the contrary, suppose that some $u_1,u_2 \in S$ has $d_G(u_1,u_2) > 2$. Note that $C_{N_G}(u_1,u_2) = \emptyset$. 

By Lemma \ref{lem:compPathLowDiam} where $u := u_1$ there exists a set $N_G(u_1) \cap A \supseteq \{a_b,a_d\}$ satisfying $b,d \in [1,k+4]$ and $\abs(b - d) > 2$. W.l.o.g. let $b < d$. 

By Lemma \ref{lem:compPathLowDiam} where $u := u_2$ there exists a set $N_G(u_2) \cap A \supseteq \{a_{p},a_{q}\}$ satisfying $p,q \in [1,k+4]$ and $\abs(p - q) > 2$. W.l.o.g. let $p < q$. Moreover, w.l.o.g. we may assume that $b < p$. Note that $b,d,p,q$ are distinct because $C_{N_G}(u_1,u_2) = \emptyset$. Notice it is possible that $b + 1 = p$, but it is not possible that $b + 1 \in \{d,q\}$. Since $b + 1 < q$, we have that $a_b \in N_G(a_q)$. We consider two cases:

\textbf{Case 1:} $p + 1 < d$. Thus, $a_p \in N_G(a_d)$. There exists a path $W = \langle a_b a_q u_2 a_p a_d \rangle$ in $G$. See the right graph in Fig. \ref{fig:ComplementPath+v} for an example of this case. However, the existence of $D_2(W) = \langle a_b u_2 a_d \rangle$ gives $d_{G_2}(a_b,a_d) = 2$, a contradiction to the assumption that $\abs(b - d) > 2$. We have proven our claim that for all pairs $v_1,v_2 \in S$ we have $d_G(v_1,v_2) \leq 2$. Since we have also shown that $G \setminus (A \cup S) = \emptyset$, it immediately follows that $\diam(G) \leq 2$, a contradiction to the hypothesis that $\diam(G) = k$.

\textbf{Case 2:} $p + 1 \geq d$. In other words, we have that $a_p,a_q \in P_2[a_{d-1},a_{k+4}]$. First, suppose that $p = b + 1$. Hence, $(b+1)+1 \geq d$. If $b + 2 \in \{d,d+1\}$, then we contradict the fact that $\abs(b - d) > 2$. So, it is necessary that $b + 2 > d + 1$, but this clearly contradicts the assumption that $b < d$. We see that $p = b + 1$ leads to contradiction. It is necessary that $p \neq b + 1$. Recalling also that $p > b$, we in fact have that $p > b + 1$ and thereby $a_b \in N_G(a_p)$. To simplify, we have shown that $a_b$ and $a_p$ are sufficiently far from one another in $P_2$ that they are adjacent in $G$.

Our assumption that $p + 1 \geq d$, combined with the fact that $\abs(p - q) > 2$ and $p < q$, implies that $q \neq d + 1$. Thus, $a_q \in N_G(a_d)$. Stated plainly, $a_q$ and $a_d$ are far enough away in $P_2$ that they are adjacent in $G$. It has become evident that there exists a path $W = \langle a_b a_p u_2 a_q a_d \rangle$ in $G$. The existence of $D_2(W) = \langle a_b u_2 a_d \rangle$ gives $d_{G_2}(a_b,a_d) = 2$, a contradiction to the assumption that $\abs(b - d) > 2$. We have proven our claim that for all pairs $v_1,v_2 \in S$ we have $d_G(v_1,v_2) \leq 2$. Since we have also shown that $G \setminus (A \cup S) = \emptyset$, it immediately follows that $\diam(G) \leq 2$, a contradiction to the hypothesis that $\diam(G) = k$.
\end{proof}

Combining our main theorem with the lower bound achieved by Jafari and Musawi \cite{Jafari}, we arrive at the following statement:

\begin{theorem}
Let $G$ be any graph with $\diam(G) = k \geq 3$. Then, either $\diam(G_2) = \infty$ or $\lceil \frac{1}{2} k \rceil \leq \diam(G_2) \leq k + 2$.
\label{thm:k+2Full}
\end{theorem}

The following can be shown with examples.
\begin{corollary}
The upper bound of the inequality expressed by Theorem \ref{thm:k+2Full} is sharp when $k = 3$ or $k > 3$ is even.
\label{cor:keven}
\end{corollary}

\begin{figure}[h]
    \centering
    ~~
        \begin{tikzpicture}[scale=1, every node/.style={scale=.5}]

        \tikzstyle{vertex}=[circle, fill=white]
        \begin{scope}\Large
        \end{scope}
        
        \tikzstyle{vertex}=[circle, fill=black]
        
        \node[vertex] (a) at (0,0) {$ $};
        \node[vertex] (b) at (0,1) {$ $};
        \node[vertex] (c) at (0,2) {$ $};
        \node[vertex] (d) at (0,3) {$ $};
        \node[vertex] (e) at (1.5,0) {$ $};
        \node[vertex] (f) at (1.5,3) {$ $};
        \node[vertex] (g) at (3,0) {$ $};
        \node[vertex] (h) at (3,3) {$ $};
        \node[vertex] (i) at (4,1.5) {$ $};

        \tikzset{decoration={snake,amplitude=.3mm,segment length=3mm, post length=0mm,pre length=0mm}}
        \draw[thick](a)--(b)--(c)--(d)--(f)--(h)--(i)--(g)--(e)--(a);
        \draw[thick](g)--(h);
        \draw[ForestGreen,bend left=45](a) to (c);
        \draw[ForestGreen,bend left=-30](a) to (g);
        \draw[ForestGreen](c) to (f);
        \draw[ForestGreen,bend left=45](a) to (c);
        \draw[ForestGreen](f) to (g);
        \draw[ForestGreen](f) to (i);
        \draw[ForestGreen](i) to (e);
        \draw[ForestGreen](b) to (e);
        \draw[ForestGreen,bend left=45](b) to (d);
        \draw[ForestGreen,bend left=30](d) to (h);
        \draw[ForestGreen](e) to (h);

        \end{tikzpicture}
        ~~
        \begin{tikzpicture}[scale=1.275, every node/.style={scale=.5},rotate=90]

        \tikzstyle{vertex}=[circle, fill=white]
        \begin{scope}\Large
            \node[vertex] (i) at (-1.55,4) {$ $};
        \end{scope}
        
        \tikzstyle{vertex}=[circle, fill=black]
        
        \node[vertex] (a) at (0,0) {$ $};
        \node[vertex] (b) at (-.5,.5) {$ $};
        \node[vertex] (c) at (.5,.5) {$ $};
        \node[vertex] (d) at (0,1) {$ $};
        \node[vertex] (e) at (0,2) {$ $};
        \node[vertex] (f) at (0,3) {$ $};
        \node[vertex] (g) at (-.5,3.5) {$ $};
        \node[vertex] (h) at (.5,3.5) {$ $};
        \node[vertex] (i) at (0,4) {$ $};

        \tikzset{decoration={snake,amplitude=.3mm,segment length=3mm, post length=0mm,pre length=0mm}}
        \draw[ForestGreen,thick](d)--(c)--(a)--(b)--(d)--(e)--(f)--(g)--(i)--(h)--(f);

        \end{tikzpicture}
        
        \includegraphics[width=0.7\linewidth]{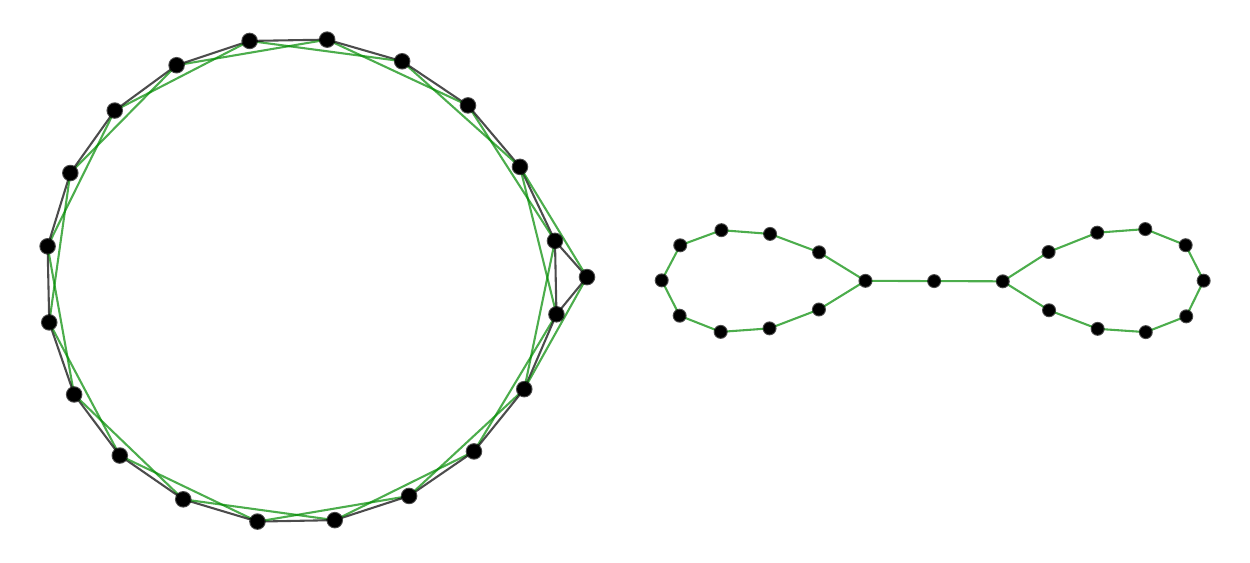}
        \caption{These graphs exemplify Corollary \ref{cor:keven} where $k = 4$ and $k = 10$. The left two graphs have $\diam(G) = 4$ and $\diam(G_2) = 6$, where black edges denote $E(G)$ and green edges denote $E(G_2)$. The right two graphs have $\diam(G) = 10$ and $\diam(G_2) = 12$.}
        \label{ref:StarCutset}

\end{figure}

\begin{proof}
An example for $k = 3$ is provided in \cite{Jafari}. There exists a family of graphs showing that sharpness holds for any even $k > 3$. See Fig. \ref{ref:StarCutset}, which shows a graph $G$ having $\diam(G) = 4$ and $\diam(G_2) = 6$. By increasing the size of the largest cycle in $G$ by $4$ such that the triangle shares exactly one edge with the largest cycle, $\diam(G)$ and $\diam(G_2)$ increase by $2$. It is easy to verify that this procedure can be repeated in order to generate any number of graphs satisfying the conclusion of the corollary.
\end{proof}


\section{Computing Higher Order 2-distance Graphs}

\begin{figure}[h]
    \centering
    \includegraphics[width=0.99\linewidth]{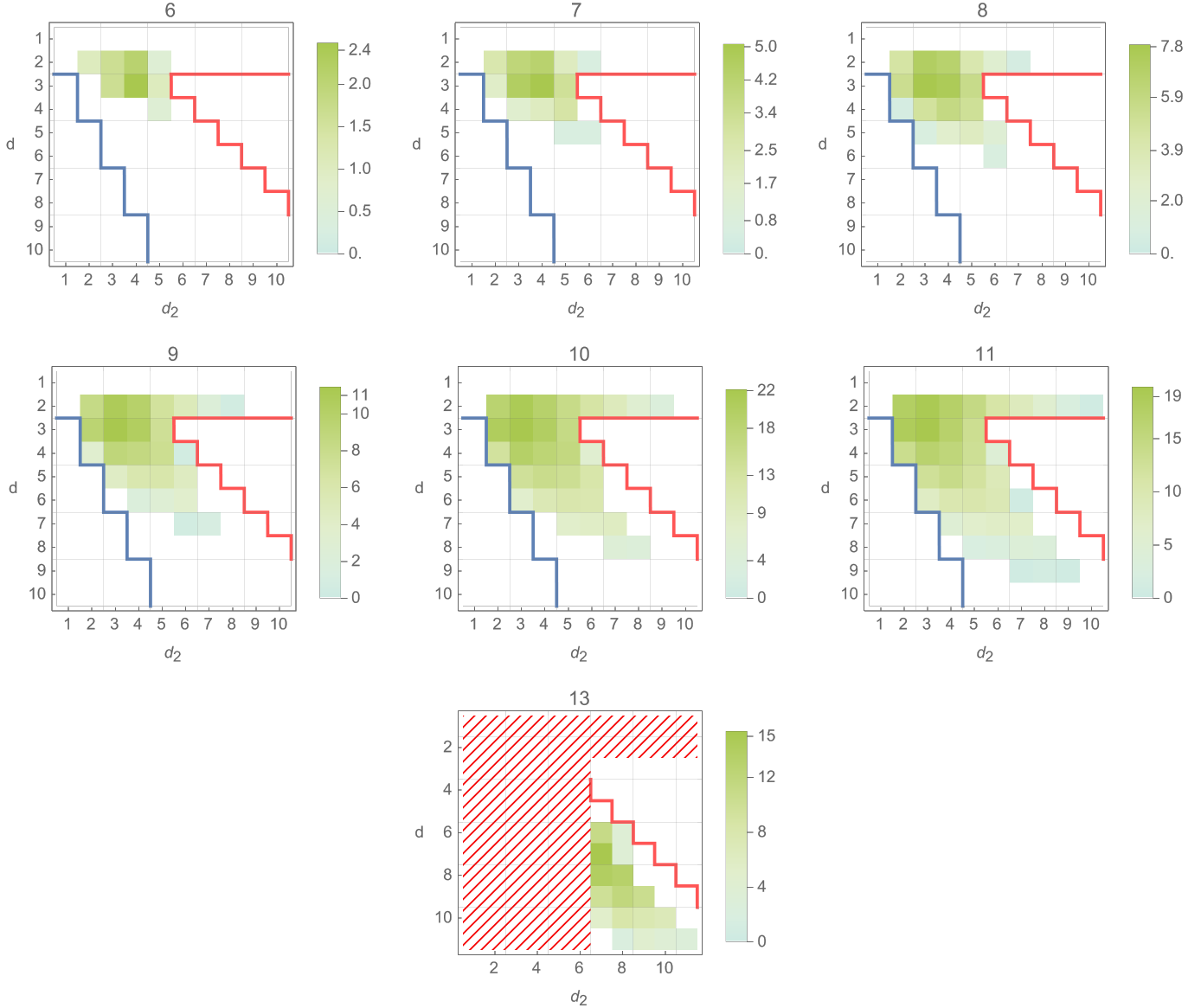}
    
    \caption{Log scale plots showing the proportion 
    of the diameters of $G$ and $G_2$ where $G$ belongs to the set of all graphs which are non-isomorphic, connected, and have fixed $|V|$. The top label corresponds to $|V|$. Axes labeled $d$ and $d_2$ denote $\diam(G)$ and $\diam(G_2)$, respectively. The blue and red boundaries express the lower and upper bounds of the inequality from Theorem \ref{thm:k+2Full}. The bottom-most plot demonstrates a large reduction in search space provided by our \textsc{sat} equations for $|V| = 13$, allowing us to find $2$-distance graphs that meet the sharp upper bound at $(d,d_2) = (6,8)$ given by Theorem \ref{thm:k+2Full}. The red region shows parameter ranges excluded by the equations. }

    \label{fig:reduced}
\end{figure}

The use of brute-force to compute all $2$-distance graphs is reliable up to $11$ vertices, with results provided in Fig. \ref{fig:reduced}. But, the run-time requirements increase exponentially beyond that. As it is apparent from the figure, graphs having $\diam(G) = 2$ and $\diam(G_2) = |V| - 1$ exist. For example, this property holds if $G$ is the complement of a path that has $|V| \geq 4$. 

Rather than using exhaustive search, $2$-distance graphs with high $|V|$ can be computed with careful use of a \textsc{sat} solver. In our implementation, we express $G = (V,E)$ and $G_2$ as adjacency matrices where $\True$ and $\False$ denote whether an edge exists. Henceforth, let $a_{i,j}$ be the $(i,j)$-th element in the adjacency matrix of $G$ and let $b_{i,j}$ be the $(i,j)$-th element in the adjacency matrix of $G_2$. By definition, 
\[
b_{i,k} = \bigvee_j (a_{i,j}\land a_{j,k} \land \lnot a_{i,k})
\] 
We fix $P_2 = \langle 0,1,2,3, \dots ,k+1,k+2 \rangle$. Then, the first set of equations for the \textsc{sat} solver fix the edges for $P_2$. That is, for all $b_{i,j} \in P_2 \textnormal{~we have~} b_{i,j} = \True$. For any $a,b \in P_2$, let $\mathscr P_{a,b,\ell}$ be the set of all $a,b$-paths of length $\ell$ in $(G_2 \setminus P_2) \cup \{a,b\}$. Then, the set of equations that forbid shortcuts that would invalidate the diametral path $P_2$ can logically be expressed as 
\[
\bigwedge_{P \in \mathscr P_{a,b,\ell}} \lnot \big(  \underset{(v,u) \in P}{\wedge} b_{v,u} \big)
\]
for all $\ell \leq |P_2|$. To eliminate the large number of isomorphic subgraphs in $\mathscr P_{a,b,\ell}$ in $H = G_2 \setminus P_2$, more edges can be fixed. Specifically, each nonisomorphic subgraph in $H$ can be assigned its own satisfiability equation. Finally, an equation to eliminate the low-diameter $k=2$ graphs is used, and is given by 
\[
\lnot \bigwedge_{i,k} (a_{i,j}\land a_{j,k})
\] 
and means that if any path $\langle a_{i,j} a_{j,k}\rangle$ does not exist (i.e. $\diam(G) \neq 2$), then this expression will evaluate to $\True$, and $\False$ otherwise.

After conversion of our logical equations to conjunctive normal form, we were able to perform an experiment restricted to $\diam(G_2) \geq 7$, $|V| = 13$, and $\diam(G) > 2$ (see Fig. \ref{fig:reduced}). As opposed to the weeks required to exhaustively enumerate graphs on $13$ vertices, we found results within one day.

\newpage

\bibliographystyle{elsarticle-num}
\bibliography{bibliography}

\end{document}